\documentclass[12pt]{article}
\usepackage{amsmath, dsfont, amssymb, amscd, caption,  color}
\usepackage{graphicx,subfigure}
\usepackage[margin=1in]{geometry}

\usepackage{times,soul}
\usepackage[square,sort,comma,numbers]{natbib}

\usepackage{setspace}
\usepackage{tikz}
\usepackage{url}
\usepackage{mathbbol}
\usepackage{tikz-cd}
\usepackage{caption}
\usepackage{amsfonts}
\usepackage{fancybox}
\usepackage{algorithmic}
\usepackage{algorithm}
\usepackage{comment}
\usepackage{color}
\usepackage[colorlinks,citecolor=blue,urlcolor=blue]{hyperref}

\usepackage{paralist}
\usepackage{tikz,enumerate}
\usepackage{multirow}
\usepackage[normalem]{ulem}

\usepackage{hyperref}
\usepackage{cleveref}

\newcommand{\bpm}{\begin{pmatrix}}
\newcommand{\epm}{\end{pmatrix}}

\graphicspath{{images/}}

\newtheorem{theorem}{Theorem}[section]
\newtheorem{lemma}[theorem]{Lemma}
\newtheorem{prop}[theorem]{Proposition}
\newtheorem{corollary}[theorem]{Corollary}

\newtheorem{definition}[theorem]{Definition}

 \numberwithin{dummy}{section}

\newenvironment{proof}[1][. ]{{\bf Proof#1}}{\hfill$\square$\vskip\baselineskip}

\begin{document}

\title{Twistor spaces on foliated manifolds}
\author{Rouzbeh Mohseni and Robert A. Wolak}

\date{}

\maketitle

\begin{abstract}
The theory of twistors on foliated manifolds is developed. We construct the twistor space of the normal bundle of a foliation. It is demonstrated that the classical constructions of the twistor theory lead to foliated objects and permit to formulate and prove foliated versions of some well-known results on holomorphic mappings. Since any orbifold can be understood as the leaf space of a suitably defined Riemannian foliation we obtain orbifold versions of the classical results as a simple consequence of the results on foliated mappings.



%\tableofcontents
\end{abstract}

%---Intro
\section{Introduction}

Let $M^{2n}$ be an even-dimensional Riemannian manifold, the twistor space $Z(M)$ is the parametrizing space for compatible almost complex structures on $M$. It is a bundle over $M$, with fiber $SO(2n)/U(n)$ and is equipped with two almost complex structures $J^\pm$, where $J^+$ can be integrable but $J^-$ is never integrable, however, it still is important as will be discussed in Theorem \ref{th5}. Moreover, in the case where $J^+$ is integrable, it is shown in \cite{AHS} that $M$ has particular properties, especially when $n=2$, which is an interesting case in physics, since the holomorphic structure of the twistor space correspond to a conformal structure of $M$. This correspondence is called the Penrose correspondence.  

In Section \ref{S2} and \ref{S3} we review some of the necessary notions on the twistor space, in particular, the construction of the twistor space $Z(M)$. Let $(M,\mathcal{F})$ be a foliated Riemannian manifold with foliation $\mathcal{F}$.   We construct the bundle of normal twistors $Z(M,\mathcal{F})$ and the bundle of tangent twistors $Z(\mathcal{F})$, hence, on a foliated Riemmanian manifold, there are  three different  twistor spaces. Twistor spaces are also important in the study of harmonic maps, and in Theorem \ref{th5} using the constructed twistor space in Section \ref{S4} we reformulate one of the classic results regarding harmonic maps. Several other results concerning harmonic maps between Riemannian foliated manifolds are formulated and proved.  

As it was remarked by A. Haefliger et al. in \cite{ha}, any orbifold can be realized as the leaf space of a foliation with compact leaves. Therefore, having the results for foliated Riemannian manifolds in Section \ref{S5}, in Section \ref{S6} we try to reformulate them for the leaf spaces, hence obtaining similar results for orbifolds.

\section{Structures on principal bundles}\label{S2}

 The arguments of this section are quite well-known but we present them to fix the notation and conventions.

Let $(M^m,g)$ be an oriented Riemannian manifold of dimension $m$ and $D$ be an oriented subbundle of $TM$ of dimension $d$. Then we denote by $SO_g(M)$ the bundle of orthonormal positive oriented frames on M. This is a principal bundle with the structure group $SO(m)$. Let $\omega$ be the Levi-Civita connection form on $SO_g(M)$. The form $\omega$ has its values in $\mathfrak{so}(m)$. There is a naturally defined Riemannian metric on $SO_g(M)$ which is defined by the pull-back of $g$ on the horizontal vectors and by the Killing form on the vertical vectors of $TSO_g(M)$.

We denote by $SO_g(D), SO_g(D^\perp)$ the principal bundles of oriented orthonormal frames of $D$ and $D^{\perp}$ with the structure groups $SO(d)$ and $SO(m-d)$, respectively. We take an orientation on $D^{\perp}$ such that elements of the oriented basis of $D$ together with that of $D^{\perp}$ give the orientation of $TM$. Then the principal fiber bundle $SO_g(D)\times_M SO_g(D^{\perp})$ has the structure group $SO(d)\times SO(m-d)$ and is canonically isomorphic to reduction of the principal bundle $SO_g(M)$ to the subbundle consisting of the frames $(v_1,..., v_d, w_1,...,w_{m-d})$ such that $(v_1,...,v_d)\in SO_g(D)_x$ and $(w_1,...,w_{m-d})\in SO_g(D^{\perp})_x$ for $x\in M$. Hence we identify $SO(D)\times_M SO(D^{\perp})$ as the subbundle of $SO(M)$.

The Lie algebra $\mathfrak{so}(m)$ may be expressed as an orthogonal sum $\mathfrak{so}(d)\oplus\mathfrak{so}(m-d)\oplus\mathfrak{m}(d,m-d)$ where the decomposition is given by

\begin{equation}{\label{1}}
    \begin{pmatrix}
    A & B\\
    -^tB & C
    \end{pmatrix}
    =
  \begin{pmatrix}
    A & 0\\
    0 & 0
    \end{pmatrix}
    +
      \begin{pmatrix}
    0 & 0\\
    0 & C
    \end{pmatrix}
    +
      \begin{pmatrix}
    0 & B\\

    -^tB & 0
    \end{pmatrix}
\end{equation}
where $A\in \mathfrak{so}(d)$, $C\in \mathfrak{so}(m-d)$ and $B\in \mathfrak{m}(d,m-d)$; Here we denote by $\mathfrak{m}(d,m-d)$ the vector space of $d\times (m-d)$ matrices with coefficients in $\mathds{R}$. Moreover, the decomposition in (\ref{1}) is orthogonal with respect to the Killing form on $\mathfrak{so}(m)$. We denote by $p_1$ the orthogonal projection $p_1 : \mathfrak{so}(m)\longrightarrow \mathfrak{so}(d)\oplus \mathfrak{so}(m-d)$.

The connection $\omega$ is usually not reducible to the subbundle $SO_g(D)\times_M SO(D^{\perp})$; This happens when both bundles $D$ and $D^{\perp}$ are totally geodesic. However, we put $\omega_1:=p_1\circ \omega$. This gives a connection on the principal bundle $SO(D)\times_M SO(D^{\perp})$. Then there exists a natural morphism of the fiber bundles $\Psi : SO(D)\times_M SO(D^{\perp})\longrightarrow SO(D)$ over the homomorphism of the Lie groups $SO(d)\times SO(m-d)\longrightarrow SO(d)$.

\begin{lemma}
{\label{lemma1}}
The connection $\omega_1$ projects via the homomorphism $\Psi$ to a metric connection on $SO(D)$.

\begin{proof}
By a straightforward  verification.
\end{proof}
\end{lemma}

We denote by $\omega_2$ the unique connection obtained in Lemma \ref{lemma1}.

The metric connection $\omega_2$ may be characterized in another way as follows. Let $\sigma=(X_1,...,X_d)$ be a local section of the bundle $SO(D)$ defined on an open subset $U$ of $M$. This means that $X_1,...,X_d$ are orthonormal vector fields on $U$ with their values in $D$. Then $\sigma^*\omega_1$ is a $1$-form with values in $\mathfrak{so}(d)$ and it has components $(\omega_2)^i_j$ defined by the following identity: 

\begin{equation}
    g(\nabla X_i,X_j)=(\omega_2)^i_j
\end{equation}
where $\nabla$ denotes the Levi-Civita connection of $g$.

The manifold $SO_g(D)$ carries a naturally defined Riemannian metric structure $g_1$ defined in the following way

\begin{align}\label{2}
    g_1:=\pi^*g_0+h_0(\omega_2,\omega_2)
\end{align}
where $h_0$ is the scaled Killing bilinear form on $\mathfrak{so}(d)$, we put $h_0(X,Y)=(1/2)trace(XY)$ for each $X,Y\in \mathfrak{so}(d)$. The metric $g_1$ is invariant with respect to the action of the group $SO(d)$ on $SO(D)$.

\section{Twistor space on a subbundle of a Riemannian manifold}\label{S3}

In this section we suppose that $(M^m,g_0)$ is a Riemannian manifold and $D$ is a subbundle of $TM$ such that the dimension of the fibers of $D$ is even and equal to $d=2q$, $(q\geq 1)$. Our discussion involves some notions of f-twistor theory and for sake of brevity some of the standard definitions are omitted, for a more detailed discussion of f-twistor spaces, f-structures, etc., see \cite{Ra}. 

Let us consider the homogeneous manifold $SO(2q)/U(q)$. This manifold is well-known in the theory of twistor spaces cf. \cite{1,2}. The structures on this manifold are essential for our paper, therefore, we recall some of the definitions and properties here. The manifold $SO(2q)/U(q)$ is a K{\"a}hlerian symmetric space. A Riemannian structure may be obtained as follows. The group $SO(2q)$ carries the standard bilinear bi-invariant Riemannian metric defined by $-(1/2)trace(XY)$ where $X$ and $Y$ are invariant vector fields identified with the anti-symmetric $d\times d$ matrices. Then the Riemannian structure on $SO(2q)/U(q)$ is obtained from the above form on $SO(2q)$ via the projection on $SO(2q)/U(q)$ in such a way that the projection is a Riemannian submersion, and we denote this metric by $h_1$. The group $SO(d)$ acts on $SO(2q)/U(q)$ on the left via isometries.

To define the almost complex structure on $SO(2q)/U(q)$ we need to decompose $\mathfrak{so}(2q)=\mathfrak{u}(q)\oplus \mathfrak{m}$ as follows
\begin{align*}
    X=(1/2)(X-J_0XJ_0)+(1/2)(X+J_0XJ_0), \,\, \forall X \in \mathfrak{so}(2q).
\end{align*}
where $J_0$ is the canonical complex structure on $R^{2q}$:  

\begin{align*}
    J_0=
    \begin{pmatrix}
    0 & -I_q\\
    I_q & 0
    \end{pmatrix}
\end{align*}
and $I_q$ is the $q\times q$ identity matrix.

The complex structure defined on $SO(2q)/U(q)$ is invariant with respect to the action of the group $SO(2q)$.

We denote by $Z(q)$ the manifold of almost complex compatible orientation preserving structures on $\mathds{R}^d$. Then there is a canonical diffeomorphism between $SO(2q)/U(q)$ and $Z(q)$. The diffeomorphism is given in the following way: if $AU(q)\in SO(2q)/U(q)$ then $AU(q)\mapsto AJ_0A^{-1}$. There is a unique Riemannian metric on $Z(q)$ such that the diffeomorphism is an isometry, also there is a natural action of the group $SO(2q)$ on $SO(2q)/U(q)$ on the left. Likewise, there is an adjoint action of $SO(2q)$ on $Z(q)$ and these actions are isometries. Furthermore, the correspondence between $SO(2q)/U(q)$ and $Z(q)$ is $SO(2q)$-equivariant. Then there is the associated bundle 

\begin{align}
    Z(D):=SO(D)\times_{SO(2q)} Z(q)
\end{align}
which is a fiber bundle over $M$ with the standard fiber $Z(q)$. Let us call it the bundle of twistors of $D$. We have the natural projection $\pi: Z(D)\longrightarrow M$ and we denote by $\mathcal{V}\longrightarrow Z(D)$ the vertical vector bundle defined by $\pi$.

The theorem below gives a natural geometric interpretation of our construction.

\begin{theorem}
The manifold $Z(D)$ parametrizes all orthogonal/compatible, orientation preserving almost complex structures on $D$. Moreover, $Z(D)$ may be also considered as a parametrization of metric f-structures on $M$ such that their image is equal to $D$ when they preserve the orientation on $D$.

\begin{proof}
Since the similar property is well-known for twistor spaces we only define the correspondence. In fact, if $[a,P]\in SO(D)\times_{SO(2q)} Z(q)$, then the corresponding almost complex structure in $D$ is given by $aPa^{-1}$. The second part of the theorem follows from the fact that in the presence of the Riemannian metric $g$ on $M$ we can extend in a unique way an almost complex structure from $D$ to a metric $f$-structure on $TM$ by posing it equal zero on $D^\bot$.
\end{proof}
\end{theorem}

Since the $SO(2q)$ action

\begin{align*}
    SO(D)\times Z(q) \times SO(2q) \rightarrow SO(D)\times Z(q)\\
    (a,P,A)\mapsto (aA,A^{-1}PA)
\end{align*}

\noindent
is an isometry then there exists the unique Riemannian metric on $Z(D)$ such that the projection

\begin{align*}
    \Phi : SO(D)\times Z(q) \rightarrow Z(D)
\end{align*}

\noindent
is a Riemannian fibration. We fix this metric on Z(D) and denote it by $\mathds{G}$.

The Riemannian structure $\mathds{G}$ can be defined in the following equivalent way. For each $a\in SO(D)$ and each $P\in Z(q)$ we consider the decomposition of tangent space
\begin{align*}
    T_{(a,P)}(SO(D)\times Z(q))= \Tilde{\mathcal{H}}_a\oplus \Tilde{\mathcal{V}}_a\oplus T_PZ(q)
\end{align*}
where $\Tilde{\mathcal{H}_a}$, $\Tilde{\mathcal{V}}_a$ are the horizontal and vertical subspaces of $T_aSO(D)$, respectively. The horizontal subspace $\Tilde{\mathcal{H}}_a$ is considered with respect to the connection $\omega_2$ on $SO(D)$.

The properties of the map $\Phi$ are essential to understand the geometric structures on $Z(D)$. We present some of them.

\begin{lemma}
Let $p\in Z(D)$ then for each $(a,P)\in p$, we have that
\begin{enumerate}
    \item $d\Phi(\Tilde{\mathcal{V}_a})=\mathcal{V}_p$
    \item $d\Phi(\Tilde{\mathcal{H}_a})$ is horizontal and independent of the choice of $(a,P)\in p$
    \item $d\Phi$ is an isomorphism when restricted to $\Tilde{\mathcal{V}_a}\oplus T_PZ(q)$.
\end{enumerate}
\begin{proof}
These properties are proved by straightforward calculations.
\end{proof}
\end{lemma}

The above lemma suggests to define the horizontal distribution $\mathcal{H}$ in $TZ(D)$ as the image of the distribution $\tilde{\mathcal{H}_a}$ via the map $d\Phi$. Moreover, we have an induced morphism of vector bundles $d\pi : \mathcal{H}\rightarrow TM$ which is an isomorphism on each fiber. Then the Riemannian structure $\mathds{G}$ on $Z(D)$ is such that: the subspaces $\mathcal{H}_p$, $\mathcal{V}_p$ are orthogonal, $\mathds{G}$ restricted to $\mathcal{H}_p$ coincides with the pull-back metric, via the projection $\pi: \mathcal{H}_p\rightarrow T_{\pi(x)}M$, and $\mathds{G}$ restricted to $\mathcal{V}_p$ is obtained from the metric of $T_PZ(q)$ via the isomorphism $d\Phi : \Tilde{\mathcal{V}}_p\rightarrow \mathcal{V}_p$. It is straightforward to prove that the metric defined in this way coincides with $\mathds{G}$.

There exist two canonically defined f-structures on $Z(D)$. We denote them by $\phi_\pm$. The f-structures $\phi_\pm$ are defined as follows: if $p=[a,P]\in Z(D)$ then for each $X\in \mathcal{H}_p$ , we consider $d\pi(X)=X_1+X_2\in D\oplus D^\bot$. Next, we put 

\begin{align*}
    \phi_\pm(X)=\pm(d\pi|_{\tilde{\mathcal{H}}_{\pi(p)}})^{-1}aPa^{-1}(X_1)
\end{align*}
the almost complex structure on $\mathcal{V}_a$ is obtained from the almost complex structure on $T_PZ(q)$ via the isomorphism $d\Phi$. Hence we have defined $\phi_\pm$ on the subspace $T_pZ(D)$. The definition of $\phi_\pm$ does not depend on the choice of the elements $(a,P)\in p$. Moreover, it is straightforward to verify that $\mathds{G}$ and $\phi_\pm$ are compatible, i.e., for each $X,Y\in T_pZ(D)$ we have that

\begin{align*}
    \mathds{G}(\phi_\pm(X),Y)+\mathds{G}(X,\phi_\pm (Y))=0.
\end{align*}

Let $\varphi: (M_1,g_1)\rightarrow (M_2,g_2)$ be a local isometry such that $d\varphi(D_1)\subset D_2$ where $D_i\subset TM_i$ are oriented subbundles of the even dimension ($i=1,2$). We assume also that $d\varphi$ preserves orientation of $D_1$ and $D_2$. Then there is the induced map $L(\varphi): SO(D_1)\rightarrow SO(D_2)$ sending, via the map $d\varphi$, the oriented frames of $D_1$ into the oriented frames of $D_2$. The map $L(\varphi)$ is $SO(d)$-equivariant and it factorizes to the map $Z(\varphi):Z(D_1)\rightarrow Z(D_2)$ which is a morphism of fiber bundles.

We have the following natural proposition which we have not found in the mathematical literature about the twistor spaces. For definition of f-holomorphic isometry, see \cite{Ra}.

\begin{prop}\label{prop1}
If $(M_1,g_1)$, $(M_2,g_2)$ are two Riemannian manifolds with even-dimensional oriented subbundles $D_1\subset TM_1$ and $D_2\subset TM_2$, then any local isometry $\varphi: M_1\rightarrow M_2$ which maps $D_1$ into $D_2$, while preserving the orientation, determines an f-holomorphic isometry $Z(\varphi): Z(D_1)\rightarrow Z(D_2)$ of the metric f-twistor spaces. We consider the same type $\pm f$-structures in the twistor spaces $Z(D_1)$ and $Z(D_2)$.

\begin{proof}
We fix $p_1\in Z(D_1)$ and $p_2 =Z(\varphi)(p_1)$. Let $(a_1,P)\in p_1$ and let $L(\varphi)(a_1)=a_2$, then $(a_2,P)\in p_2$. Since $\varphi$ is a local isometry, it sends the horizontal subspace $\Tilde{\mathcal{H}}_{a_{1}}$ of $T_{a_{1}} SO(D_1)$ onto the horizontal subspace $\Tilde{\mathcal{H}}_{a_{2}}$ in $T_{a_{2}} SO(D_2)$. Then we have the map

\begin{align}
    dL(\varphi)\times id: \Tilde{\mathcal{H}}_{a_1}\oplus T_PZ(q)\rightarrow \Tilde{\mathcal{H}}_{a_2} \oplus T_PZ(q)
\end{align}

\noindent
which determines, via the isomorphisms $d\Phi_i: \Tilde{\mathcal{H}}_{a}\oplus T_PZ(q)\rightarrow T_{p_i}Z(D_i)$ $(i=1,2)$, the map $dZ(\varphi): T_{p_{1}}Z(D_1)\rightarrow T_{p_{2}}Z(D_2)$. It is enough to prove that the map $dL(\varphi)\times id$ restricted to $\Tilde{\mathcal{H}}_{a_1}\times T_PZ(q)$ is an f-holomorphic isometry. It is clear that the map $id: T_PZ(q)\rightarrow T_PZ(q)$ is an f-holomorphic isometry.

Then we observe that the restriction of $d_aL(\varphi)$ to $\Tilde{\mathcal{H}}_{a_1}$ equals to $(d\pi_2|_{\Tilde{\mathcal{H}}_{a_2}})d\varphi\circ d\pi_1$. Thus $d_{a_1}L(\varphi)$ is an isometry as a composition of isometries.

If $X\in \Tilde{\mathcal{H}}_{a_1}$ then $d\pi_1(X)=av+X_0$, where $v\in \mathds{R}^d$ and $X_0\in D^\bot$. Then 

\begin{equation*}
\begin{split}
     dL(\varphi)\tilde{\phi}_\pm(X) & =\pm(d\pi)^{-1}(d\varphi a_1P(v))\\
   & =\pm(d\pi_2)^{-1}(a_2\overbrace{a_2^{-1}d\varphi a_1}^{=id}P(v))\\
   & =\tilde{\phi}_\pm(d\pi_2)^{-1}(a_2v+d\varphi (X_0))\\
   & =\tilde{\phi}_\pm dL(\varphi)(X).
   \end{split}
\end{equation*}
Hence $L(\varphi)$ is f-holomorphic and then $Z(\varphi)$ is f-holomorphic too.
\end{proof}
\end{prop}

Let $(M^{2n},g)$ be a $2n$-dimensional oriented Riemannian manifold and $SO_g(M)$ be the bundle of orthonormal positive frames on $M$. The twistor bundle $Z(M)$ is a bundle over $M$ with fibre $SO(2n)/U(n)$ and structure grouop $SO(2n)$ and is the associated bundle to $SO_g(M)$, we invite the interested reader to see \cite{1} for more details. Now using this and Proposition \ref{prop1} we have the following corollary which can be proved in similar steps.  

\begin{corollary}\label{C1}
Let $(M^{2n}_1,g_1)$ and $(M^{2n}_2,g_2)$ be two oriented Riemannian manifolds and $\varphi: M_1\rightarrow M_2$ be an orientation preserving local isometry, then the map $Z(\varphi): Z(M_1)\rightarrow Z(M_2)$ is a local isometry which is holomorphic with respect to the corresponding almost complex structures.
\end{corollary}

By taking into consideration results of \cite{1} we also get the following corollary.

\begin{corollary}\label{C2}
Let $(M^{2n}_1,g_1)$ and $(M^{2n}_2,g_2)$ be two oriented Riemannian manifolds and $\varphi: M_1\rightarrow M_2$, be an orientation preserving local diffeomorphism, then there exists a fiber morphism $\Psi: SO(M_1)\rightarrow SO(M_2)$ which is holomorphic if and only if $\varphi$ is conformal.

\begin{proof}
We consider on $M_1$ , the metric $\varphi^*g_2$. Then $L(\varphi): Z(TM_1,\varphi^*g_2)\rightarrow Z(TM_2,g_2)$ is a holomorphic isomorphism with respect to the induced Riemannian and complex structures. Then there exists an isomorphism of the bundles $\hat{\Phi}:Z(TM_1,g_1)\rightarrow Z(TM_1,\varphi^*g_2)$, cf. Proposition 2.2.6 of \cite{1} . Hence, by considering the composition $L(\varphi)\circ \hat{\Phi}$ and applying Proposition 2.2.7 of \cite{1} we get our corollary.

\end{proof}
\end{corollary}

Let $\rho: (M,g)\rightarrow (B,h)$ be a Riemannian submersion and $B$ be even-dimensional and oriented. Let $D=(kerd\rho)^\bot$ be the horizontal distribution with its natural orientation taken, via $\rho$, from $B$. Then for each $x\in M$ and each $P_x\in Z(D)_x$ we put
\begin{align*}
    \Xi (P_x)=d_x\rho\circ P_x\circ (d_x\rho \vert_D)^{-1}.
\end{align*}
We equip both $Z(D)$ and $Z(B)$ with the Riemannian metrics and the almost complex structures coming from the twistor construction. We have the following lemma which considers Riemannian foliations and twistor spaces.

\begin{lemma} $\Xi$ is a Riemannian fibration and is f-holomorphic with respect to the same type of metric f-structures. Moreover, the following diagram
\begin{displaymath}
  \begin{tikzcd}
   Z(D) \arrow[r, "\Xi "] \arrow[d]
   & Z(B) \arrow[d]\\
   M \arrow[r, "\rho"]
   & B
 \end{tikzcd}
 \end{displaymath}
 commutes.
 
 \begin{proof}
 Let $\omega_B$ be the Levi-Civita connection form on $SO(B)$. Then $\rho^*\omega_B$ is the connection on $SO(D)$ obtained in Lemma 2.1. We have the induced map $L(\rho): SO(D)\rightarrow SO(B)$ such that $L(\rho)(v_1,...,v_d)= (d\rho(v_1),...,d\rho(v_d))$, where $(v_1,...,v_d)\in SO(D)$. The map $L(\rho)$ is a diffeomorphism on the fibers of these principal bundles. Since $dL(\rho)$ sends the horizontal subbundle of $SO(D)$ onto the horizontal subbundle of $SO(B)$. We denote by $\pi_1:SO(D)\rightarrow M$ and $\pi_2: SO(B)\rightarrow B$ the canonical projections. Let $a\in SO(D)_x$ and $b=L(\rho)(a)\in SO(B)_{\rho(x)}$. We denote by $\Tilde{\mathcal{H}}_a$ the horizontal space in $T_aSO(D)$ and by $\hat{\mathcal{H}}_b$ the horizontal space in $T_bSO(D)$. Then the maps $d\pi_1:\tilde{\mathcal{H}}_a\rightarrow T_xM$, $d\pi_2:\hat{\mathcal{H}}_b\rightarrow T_{\rho(x)}B$ are isomorphisms. Then for each $P\in Z(q)$, we have the following commutative diagram
 \begin{displaymath}
  \begin{tikzcd}
   \Tilde{\mathcal{H}}^\prime_a\oplus\tilde{\mathcal{H}}^{''}_a \arrow[r, "dL(\rho) "] \arrow[d, "d\pi_1"]
   & \hat{\mathcal{H}}_b \arrow[d, "d\pi_2"]\\
   (\mathit{ker} \, d_x\rho)\oplus D \arrow[r, "d_x\rho"]
   & T_{\rho(x)}B
 \end{tikzcd}
 \end{displaymath}
 Hence, it follows that $dL(\rho)$ sends $\tilde{\mathcal{H}}^{''}_a$ isomorphically onto $\hat{H}_b$ and has kernel equal to $\tilde{H}^\prime_a$. Therefore, for each $P\in Z(q)$ we have that
 
 \begin{align}\label{3}
     dL(\rho) \times id:  \Tilde{\mathcal{H}}^\prime_a\oplus\tilde{\mathcal{H}}^{''}_a\oplus T_PZ(q)\rightarrow \hat{H}_b\oplus T_PZ(q)
 \end{align}
 is an isometry when restricted to $\tilde{\mathcal{H}}^{''}_a\oplus T_PZ(q)$, has kernel $\tilde{\mathcal{H}}^\prime_a$ and is f-holomorphic if we equip the spaces in (\ref{3}) with their metric and holomorphic structure. Then after the projection on the twistor spaces $T_{[a,P]}Z(D)\rightarrow T_{[b,P]}Z(B)$, we get the desired result.
 \end{proof}
\end{lemma}

\section{Twistor spaces on foliated manifolds}\label{S4}

Let $(M,\mathcal{F})$ be a Riemannian foliation. Then it is defined by a cocycle $\mathcal{U}=\{U_i,f_i,k_{ij}\}_{i,j\in I}$ that is modelled on a Riemannian manifold $(N,\Bar{g})$ such that 
\begin{enumerate}
    \item $f_i:U_i\rightarrow N$ is a submersion with connected fibers;
    \item $k_{ij}: f_j(U_i\cap U_j)\rightarrow f_i(U_i\cap U_j)$ are local isometries of $(N,\Bar{g})$;
    \item $f_i=k_{ij}f_j$ on $U_i\cap U_j$.
\end{enumerate}

Let $Q=TM/T\mathcal{F}$ be the normal bundle of the foliation $\mathcal{F}$. The natural action of the Lie algebra sheaf $({\mathcal X},{\mathcal F})$ of local vectors tangent to the leaves of the foliation $\mathcal F$  defines the foliation ${\mathcal F}_Q$ of the manifold $Q$ which is of the same dimension as $\mathcal F$. Its leaves are covering spaces of the corresponding leaves of $\mathcal F$. Therefore the tensor bundle $\otimes^2 Q^*$ is a foliated bundle with the induced foliation ${\mathcal F}_Q^2$ of the same dimension as $\mathcal F$.  Over the open set $U_i$ the normal bundle is isomorphic to the pull-back $f^*_iTN$. This isomorphism sends traces of leaves of ${\mathcal F}_Q$ to the fibres of the submersion $df_i\vert Q$. Moreover, it induces an isomorphim of the foliated  bundles 
$\otimes^2Q^*_{\vert U_i}$ and $f_i^*\otimes^2TN^* $. On the manifold $M$ we have a bundle-like metric $g$ which induces a metric tensor $g_Q$ on $Q$ which is a foliated section of $Q$, i.e., a foliated mapping 
$g_Q \colon (M,{\mathcal F}) \rightarrow (Q^* \otimes Q^*,{\mathcal F}^2_Q)$. Under the same isomorphisms the tensor $g_Q$ corresponds to $f_{i}^{*}\bar{g}$. Summing up, we can say that a  Riemannian metric $g_Q$ in the normal bundle $Q$ is foliated iff for any foliated local sections $X,Y$ of $Q$ the function $g(X,Y)$ is basic. 

At any point $x\in M$ we define the normal twistors as the set
\begin{align}\label{Z}
    Z(M,\mathcal{F})_x=\{J^Q_x:Q_x\rightarrow Q_x, (J^Q)^2=-id_Q, \;\;\; g_Q(J^Q,J^Q) =g_Q\}
\end{align}
and the bundle of normal twistors $Z(M,\mathcal{F})$ as the union of $Z(M,\mathcal{F})_x$ over all $x\in M$.

The cocycle definition of the foliation permits us to relate normal twistors on the foliated manifold $(M,\mathcal{F})$ to the twistors of the model manifold $N$, which will be called transverse twistors. In fact, at any point $x\in M$, the differential $d_xf_i$ defines a linear isomorphism $d_xf_i|_{Q_x}: Q_x \rightarrow T_{f_i(x)}N$, which in turn defines an isomorphism of the twistor spaces $Z(M,\mathcal{F})_x$ and $Z(N)_{f_i(x)}$. Therefore any submersion $f_i$ defines a submersion $Z(f_i)$ on the level of twistors 
\begin{align*}
    Z(f_i): Z(M,\mathcal{F})|U_i\rightarrow Z(N).
\end{align*}
Hence the bundle $Z_i(M,\mathcal{F}):=Z(M,\mathcal{F})|U_i$ is isomorphic to $f^*_iZ(N)$. Moreover, the cocycle $Z(\mathcal{U})=\{Z_i(M,\mathcal{F}),Z(f_i),Z(k_{ij})\}$ defines a foliation $\mathcal{F}_Z$ on the manifold $Z(M,\mathcal{F})$, which is modelled on $Z(N)$ and the leaves of $\mathcal{F}_Z$ are covering spaces of leaves of $\mathcal{F}$.

The manifold $Z(N)$ carries the natural almost Hermitian structure $(h_N, J_N)$ and transformations $Z(k_{ij})$, according to Corollary \ref{C1}, are holomorphic isometries of the almost Hermitian structure of $Z(N)$. The normal bundle $Q_Z$ of the foliation ${\mathcal F}_Z$, locally, is isomorphic to the pullback $Z(f_i)^*Z(N)$ as a foliated manifold. Therefore the bundle $Q_Z$ admits a natural Hermitian metric $(h_Z,J_Z)$ which, in fact, is foliated and the submersions $Z(f_i)$ are transversely  Hermitian for any extension of $h_Z$ to a Riemannian metric on $M$. Moreover, the mappings $Z(f_i)$ are almost Hermitian submersions. Thus the foliation $\mathcal{F}_Z$ on the manifold $Z(M,\mathcal{F})$ is transversely  almost Hermitian for the canonical Riemannian metric and the metric compatible transverse almost complex structure $J^Q$. The normal Hermitian metric can be constructed directly without involving the transverse manifold and its twistor bundle. 

As in \cite{Ba} one can construct a twistor bundle of a vector bundle. The normal bundle is a foliated vector bundle over $M$ with foliation $\mathcal{F}_Q$ described earlier and after performing the twistor construction a foliated twistor bundle will be obtained which we denote by $Z(Q,\mathcal{F}_Q)$, moreover, using the model manifold $N$, the twistor space of $Q$ can be constructed that is isomorphic to the pull-back bundle of $Z(N)$ and has a cocycle which takes values in $Z(\mathcal{U})$. The important fact is that the twistor bundle of the normal bundle constructed in either way is isomorphic to the other.

The fact that the lifted foliation $\mathcal{F}_Z$ on $Z(M,\mathcal{F})$ is given by the lifted cocycle $Z(\mathcal{U})$ can be used to prove some interesting properties of Riemannian foliations using the well-developed theory of twistor spaces. Let $\Omega:\bigwedge^2TM\longrightarrow \bigwedge^2TM$ be the curvature operator \cite{1} of $(M,\mathcal{F})$ defined as
\begin{align*}
    \langle \Omega(X\wedge Y),Z\wedge W\rangle=\langle R(X,Y)Z,W\rangle
\end{align*}
where $X,Y,Z,W\in \Gamma(TM)$ are tangent vector fields over manifold $M$ and $R$ is the Riemann curvature tensor. Having $Q$ the normal bundle of the foliation , one can define the basic curvature operator $\Omega_b:\bigwedge^2Q\longrightarrow \bigwedge^2Q$ similarly using the transversal curvature tensor $R^Q$ \cite{Ro} defined as the curvature of the normal Levi-Civita connection $\nabla^Q$ given by the formula $2g_Q(\nabla^Q_X Y,Z)=Xg_Q(Y,Z)+Yg_Q(Z,X)-Zg_Q(X,Y)+g_Q([X,Y],Z)-g_Q([Y,Z],X)+g_Q([Z,X],Y)$, with $X,Y,Z\in \Gamma(Q)$, instead of the Riemann curvature tensor $R$.

\begin{theorem}\label{TH2}
Let $\mathcal{F}$ be a Riemannian foliation of codimension $2q\geq 4$ on a compact connected  manifold $M$. If $\mathcal{F}$ admits a compact transverse manifold $S$ and the basic curvature operator $\Omega_b$ is the identity, then the foliation $\mathcal{F}$ is developable.

\begin{proof}
Since the foliation is Riemannian, the closures of leaves are compact submanifolds. Therefore the saturation of the transverse manifold $S$ is both open and compact, thus the whole manifold, which shows that  $S$ is a complete transverse manifold. Hence foliation $\mathcal{F}$ can be modelled on $(S,\Bar{g})$, where $\Bar{g}$ is the induced Riemannian metric, as a result, the lifted foliation is modelled on $Z(S)$. The operator $\Omega_b$ induces the operator $\Bar{\Omega}$ for the metric $\Bar{g}$. Our assumption assures that $\Bar{\Omega} = id$. Then Proposition 3.2.1 of \cite{1} ensures that the manifold $Z(S)$ with the canonical structure is a K{\"a}hler manifold, and therefore, according to Theorem 4.3.2 of \cite{1} , the manifold $S$ is conformally equivalent to $S^{2q}$ with the standard Riemannian structure. Finally, $\mathcal{F}$ is a transversely conformal foliation modelled on $S^{2q}$ with the standard Riemannian structure, and so it is a developable foliation.
\end{proof}
\end{theorem}

\section{Harmonic mappings}\label{S5}

In this section we prove some foliated versions of the results known in complex differential geometry. As described in \cite{w1}, for consideration of geometric structures on foliated manifolds, the two approaches named foliated and transverse are in one-to-one correspondence, i.e., the foliated objects have holonomy invariant counterparts on the transverse manifold. 

Let $(M,g,\mathcal{F})$ be a foliated Riemannian manifold with foliation $\mathcal{F}$ of codimension $m=2q$ and the normal bundle $Q=TM/T\mathcal{F}$, and the projection map $\pi:TM\longrightarrow Q$. Using the decomposition of the tangent bundle $TM=T\mathcal{F}\oplus T\mathcal{F}^\perp$, and identifying $T\mathcal{F}^\perp$ with $Q$, a metric $g_Q$ is induced on $Q$. Foliation $\mathcal{F}$ is called a transversely almost Hermitian foliation, if it is equipped with a foliated metric $g_Q$ and a transverse almost complex structure $(J^Q)^2=-Id_Q$ such that 

\begin{equation*}
    g_Q(JX, JY)=g_Q(X,Y),\,\, \forall X,Y\in \Gamma(Q).
\end{equation*}
\noindent
Furthermore, having $\pi(X^\prime)=X$ and $\pi(Y^\prime)=Y$ for each $X^\prime,Y^\prime\in \Gamma(TM)$, a normal K{\"a}hler 2-form $w$ and its dual Fundamental 2-vector $F$ can be defined as follows

\begin{align}
    w(X^\prime,Y^\prime)=g_Q(JX,Y),\\
    g_Q(F, X\wedge Y)=g_Q(JX,Y)=-g_Q(X,JY).\label{f1}
\end{align}
where $F\in\Gamma(M,\bigwedge^2 Q)$. If the 2-form $\omega$ is a foliated section of $\bigwedge^2Q^*$ or equivalently defines  a basic 2-form on $(M,{\mathcal F})$, then the 2-vector $F$ is a foliated section of $(\bigwedge^2 Q, {\mathcal F}_Q)$ and it corresponds to a $2$-vector on the transverse manifold.

Harmonic and holomorphic maps between Riemannian manifolds have been a rich subject for research. Now, we want to formulate the counterparts of a few of those results in the case of foliated Riemannian manifolds. In order to do so, we need to formulate some foliated versions of the notions already well-known for the Riemannian manifolds. The theory of transversely  harmonic mappings has been developed by the second author in cooperation with the late Jerzy Konderak in \cite{kw1, kw} and foliated mappings preserving transverse (almost) complex structures have been studied by many authors, e.g., \cite{dip}.

\begin{definition}\label{def1}
Let $(M,\mathcal{F},J^Q)$ and $(M^\prime,\mathcal{F}^\prime,J^{Q^{\prime}})$ be two foliated Riemannian manifolds having transverse almost complex structures $J^Q$ and $J^{Q^{\prime}}$ respectively. The smooth foliated (i.e., leaf-preserving) map $\phi:(M,\mathcal{F})\longrightarrow (M^\prime,\mathcal{F}^\prime)$ is transversely  holomorphic if 
\begin{align*}
    \phi^{T}_{*}\circ J^Q=J^{Q^{\prime}}\circ \phi^{T}_{*}.
\end{align*}
where $\phi^{T}_{*}=\pi^\prime\circ \phi_{*}$ , is defined as the composition of the push-forward map $\phi_{*}$ and the projection $\pi^\prime:TM^\prime\longrightarrow Q^\prime$ and it is a foliated section of the bundle $Q^*\otimes \phi^{-1}Q^\prime\longrightarrow M$.
\end{definition}

Let $\{E_a\}_{a=1,\ldots,m}$ be a local orthonormal basis frame of the normal bundle Q and the transverse tension field of $\phi$ be defined as 

\begin{align*}
    \tau_b(\phi)=tr_Q\tilde{\nabla}\phi^{T}_{*}
\end{align*}
where $\tilde{\nabla}$ is the connection on the bundle $Q^*\otimes \phi^{-1}Q^\prime$ induced by the normal Levi-Civita connection. $\tau(\phi)$ is a foliated section of the pull-back bundle $\phi^{-1}Q^\prime$. 

\begin{definition}
Let $\phi:(M,g,\mathcal{F},J^{Q}) \longrightarrow (M^\prime,g^\prime,\mathcal{F}^\prime,J^{Q^{\prime}})$ be a foliated map between two transversely  almost Hermitian foliated manifolds. It is transversely  harmonic if the transverse tension field $\tau_b(\phi)$ vanishes.
\end{definition}

In \cite{kw1} it has been demonstrated that this definition is equivalent to the condition that for any cocycles defining the foliations the induced mapping of the transverse manifolds is harmonic.   

Let $Q^\mathds{C}=Q\otimes \mathds{C}$ be the complexified normal bundle that is naturally a foliated manifold and using the transverse almost complex structure $J^Q$ from $(M,\mathcal{F})$, which maps foliated sections of the normal bundle to foliated sections of the normal bundle, we have the decomposition $Q^\mathds{C}=Q^{1,0}\oplus Q^{0,1}$ that is also compatible with the foliation and
\begin{align*}
    Q^{1,0}=\{X\in \Gamma Q^\mathds{C}|J^QX=iX\}\\
    Q^{0,1}=\{X\in \Gamma Q^\mathds{C}|J^QX=-iX\}
\end{align*}

\noindent
The $k$th exterior power is denoted as
\begin{align}
    (\wedge^k Q)^\mathds{C}=\bigoplus_{p+q=k} Q^{p,q}
\end{align}
where $Q^{p,q}\cong \wedge^p(Q^{1,0})\otimes \wedge^q(Q^{0,1})$, moreover, we have the natural isomorphisms $(Q^{p,q})^*\cong Q^{\overline{p,q}}\cong Q^{q,p}$. Let $\alpha$ be the second fundamental form of the subspace $Q^{1,0}$ in $Q^\mathds{C}$, it can be shown that it belongs to the subspace $Q^{0,2}$ and is related to $\nabla_X F$, where $\nabla$ is the Levi-Civita connection and $X\in Q^\mathds{C}$, we have 
\begin{align}
    Q^{0,2}\otimes (Q^*)^\mathds{C}\cong Q^{0,2}\otimes (Q^{0,1}\oplus Q^{1,0})\\\nonumber
    \cong (Q^{0,1}\otimes Q^{0,2})\oplus Q^{1,2}
\end{align}
The covariant derivative of the fundamental 2-vector $F$ (\ref{f1}) can be written as $\nabla F=D_1F+D_2F$ (cf. Section 1\cite{sa}), where

\begin{align*}
    D_1F\in (Q^{1,0}\otimes Q^{2,0}) \oplus (Q^{0,1}\otimes Q^{0,2}),\\
    D_2F\in (Q^{1,2}\oplus Q^{2,1}).
\end{align*}

An almost complex structure $J$ on the normal bundle $Q$ is said to be integrable if its Nijenhuis tensor $N_J$ vanishes, i.e. for any sections $X,Y\in \Gamma Q\,$, $N_J(X,Y) = [J(X),J(Y)] - J([X,J(Y)] - J([J(X),Y]) - [X,Y] =0$ where $[X,Y] = \pi [X^\prime,Y^\prime]$, $X^\prime,Y^\prime\in \Gamma(TM)$.
The transverse almost complex structure $J^Q$ is foliated and correspondingly it induces an almost complex structure $J^N$ on the transverse manifold $N$. The foliated structure $J^Q$ is integrable iff the induced almost complex structure $J^N$ is integrable. 

\begin{lemma}\label{Lemma 1.2}
The following holds
    \begin{itemize}
        \item $D_1F=0 \Longleftrightarrow \nabla_X(Q^{1,0})\subset Q^{1,0},\,\, \forall X\in \Gamma(M,Q^{1,0})$
        \item $D_2F=0 \Longleftrightarrow \nabla_X(Q^{1,0})\subset Q^{1,0},\,\, \forall X\in \Gamma(M,Q^{0,1})$
    \end{itemize}

\end{lemma}

\begin{definition}
A transversely  almost Hermitian foliation is called transversely  $(1,2)$-symplectic if $D_2F=0$ and transversely  cosymplectic if $tr(D_2F)=0$.
\end{definition}

Trivially transversely $(1,2)$-symplectic condition implies being transversely cosymplectic. Now we can state the following theorem, which has its roots in the work done by A. Lichnerowicz \cite{Li}.

\begin{theorem}\label{th 3}
Let $(M, g, \mathcal{F}, J^Q)$ and $(M^\prime, g^\prime, \mathcal{F}^\prime, J^{Q^{\prime}})$ be transversely  almost Hermitian foliated manifolds, that are respectively transversely  cosymplectic and transversely  $(1,2)$-symplectic. Then any transversely  holomorphic map $\phi:M\longrightarrow M^\prime$ is transversely  harmonic.

\begin{proof}
Let $\{\alpha_j\}$ be a local frame for $Q^{1,0}$ and as $M$ is cosymplectic we have $tr(D_2F^M)=0$. The transverse tension field is 
\begin{align*}
    \tau_b(\phi)=tr_Q\tilde{\nabla}\phi^{T}_{*}=\sum_{j} (\tilde{\nabla}_{\Bar{\alpha_j}}\phi^{T}_{*})(\alpha_j)=\sum_j \nabla^{M^\prime}_{\phi^{T}_{*}\Bar{\alpha_j}}(\phi^{T}_{*}\alpha_j)-\phi^{T}_{*}(\nabla^M_{\Bar{\alpha_j}}\alpha_j)
\end{align*}
which using the fact that $D_2F^{M^\prime}=0$ and Lemma \ref{Lemma 1.2} , it can be shown that it belongs to ${Q^{\prime}}^{1,0}$, however, since $\tau(\phi)$ is real it must vanish.
\end{proof}
\end{theorem}

Let us recall that on a foliated Riemannian manifold $(M^n,\mathcal{F})$ with a bundle-like metric and foliation of codimension $q$, there exist two Hodge star operators 

\begin{enumerate}
    \item $*:\Omega^k(M)\longrightarrow \Omega^{n-k}(M)$
    \item $\star: \Omega^k(M,\mathcal{F}) \longrightarrow \Omega^{q-k}(M,\mathcal{F}) $
\end{enumerate}
where $\Omega^k(M)$ and $\Omega^k(M,\mathcal{F})$ denote the vector spaces of $k$-forms and basic $k$-forms respectively. The Hodge star $\star$ acts on basic forms and is determined by the transverse part of the Riemannian metric. Now let us denote by $d_B$ the restriction of the differential operator $d$ to basic forms, and by $\kappa$ the mean curvature form of the foliation $\mathcal{F}$. We then have two codifferential operators as follows

\begin{enumerate}
    \item $\delta=*d*:\Omega^k(M)\longrightarrow \Omega^{k-1}(M)$
    \item $\delta_B=(d_B-\kappa\wedge)^\star:\Omega^k(M,\mathcal{F})\longrightarrow \Omega^{k-1}(M,\mathcal{F})$
\end{enumerate}
where $(d_B-\kappa\wedge)^\star\beta=(-1)^{q(k+1)+1}\star(d_B-\kappa\wedge)\star\beta$, for any $\beta\in \Omega^k(M,\mathcal{F})$. As a remark, if we have $\delta \beta=0$ for a basic form $\beta$ then we get $\delta_B=0$, for more details see \cite{kw}.

Similarly, one can instead of the fundamental $2$-vector $F$ work with the K{\"a}hler form $w$, hence, the condition for cosymplecticity will be given as $\delta w=0$. Also, the condition for being $(1,2)$-symplectic will be $(dw)^{1,2}=0$ , i.e., being $(1,2)$-closed, hence the name $(1,2)$-symplectic. Theorem \ref{th 3} has also been proved in \cite{kw}, while using these notions.

Let $\nabla$ be a connection in the normal bundle $Q$. Denote by $\omega$ its connection form on the orthonormal frame bundle $O(Q),$ which is an $O(2q)$ principal fibre bundle, where $2q$ is the codimension of the foliation. $\Gamma = ker \omega$ is an $O(2q)$-invariant  horizontal subbundle. The vector bundle $Q$ can be understood as an associated bundle to $O(Q)$ with the standard fibre $R^{2q}$. As the horizontal subbundle is  $O(2q)$-invariant,  it induces a horizontal subbundle $H$ in $Q$. The total  spaces of both bundles carry the natural foliations ${\mathcal F}^L$ and ${\mathcal F}_Q$, respectively, whose leaves are covering spaces of the initial  foliation $\mathcal F$. The identification of $Q$ with the associated bundle of $O(Q)$ identifies the foliations as well. If $(M,{\mathcal F})$ is a foliated Riemannian manifold, its Levi-Civita connection induces a foliated connection in $O(Q)$. Its connection form $\omega$ is basic on the foliated manifold $(O(Q),{\mathcal F}_L)$ which means that for any vector $v\in T{\mathcal F}_L$, $\omega (v) =0$, hence the subbundle  $T{\mathcal F}_L \subset \Gamma$. Therefore for the normal bundle $Q$ we have $T{\mathcal F}_L \subset H$.

Let us denote the normal twistor bundle by $W:=Z(M,\mathcal{F})$. It is the associated bundle of $O(Q)$ with the standard fiber $O(2q)/U(q)$. The above considerations are also valid in this case. Then we have the splitting 
\begin{align*}
    TW=H^{(\omega)}\oplus V
\end{align*}
where $H^{(\omega)}$ is the horizontal subbundle of the tangent bundle and depends on the choice of $\omega$. For the Levi-Civita connection we denote the horizontal subbundle by $H$. The horizontal subbundle $H$ splits into $H_{\mathcal F} \oplus H_Q$ which corresponds to the splitting of the tangent bundle $TM$. As the connection is bundle-like,  $H_{\mathcal F}$ is just the subbundle $T{\mathcal F}_Q$. 

Let $x\in M$, we can associate to any $J^Q\in W_x$ , its fundamental 2-vector $F\in\bigwedge^2Q_x$, which gives 
\begin{align}
    i:W\longrightarrow \bigwedge\nolimits^2Q.\label{i}
\end{align}
Recall the projection $\pi:W\longrightarrow M$ and also consider the pull-back bundle $\pi^{-1}Q$ , its fiber at $J^Q\in W_x$ is isomorphic to $Q_x$ and using $J^Q$ it is decomposed as

\begin{align*}
    (\pi^{-1}Q)^\mathds{C}= {Q^\prime}^{1,0}\oplus {Q^\prime}^{0,1}
\end{align*}
where the decomposition is compatible with the foliation since $J^Q$ is foliated, ${Q^\prime}^{1,0}$ and ${Q^\prime}^{0,1}$ are subbundles over $W$ . Moreover, on the tangent bundle of the normal twistor space $TW=T\mathcal{F^\prime}\oplus \tilde{Q}$, where $T\mathcal{F^\prime}$ is the tangent bundle to leaves of the foliation $\mathcal{F^\prime}$ of the normal twistor space and $\tilde{Q}$ is the normal bundle. Hence there exists similar decomposition on the horizontal subbundle $H$ and let us denote by $H_{\tilde{Q}}$ the projected part to the normal bundle, using the inverse of the isomorphism $H_{\tilde{Q}}\longrightarrow \pi^{-1}Q$ we have 
\begin{align}
    H_{\tilde{Q}}^\mathds{C}=({Q^\prime}^{1,0})^h\oplus ({Q^\prime}^{0,1})^h
\end{align}
which gives us an almost complex structure $J^h$ on $H_{\tilde{Q}}$. Furthermore, consider the following decomposition which can be obtained by varying $J^Q\in W$

\begin{align*}
    \left(\pi^{-1}\bigwedge\nolimits^2Q\right)^\mathds{C}={Q^\prime}^{1,1}\oplus ({Q^\prime}^{2,0}\oplus {Q^\prime}^{0,2})
\end{align*}
moreover, it can be shown that for the subbundle $V_{\tilde{Q}}$ which is the projection of the vertical subspace $V$ to the normal subspace $\tilde{Q}$, we have

\begin{align*}
    V_{\tilde{Q}}^\mathds{C}=({Q^\prime}^{2,0})^v \oplus ({Q^\prime}^{0,2})^v
\end{align*}
which defines an almost complex structure $J^v$, such that 
\begin{align*}
    J^vX=iX,\,\, \forall X\in \Gamma ({Q^\prime}^{2,0})^v,\\
    J^vX=-iX,\,\, \forall X\in \Gamma ({Q^\prime}^{0,2})^v.\\
    \end{align*}
as a result on the normal bundle of the canonical foliation of $W$ denoted by $Q(W)$, there exist two distinct almost complex structures $J^+=J^h\oplus J^v$ and $J^-=J^h\oplus (-J^v)$. For the transverse vectors belonging to $(1,0)$- class we have the following
\begin{align*}
    Q^{1,0}(W,J^+)=({Q^\prime}^{1,0})^h \oplus ({Q^\prime}^{2,0})^v\\
    Q^{1,0}(W,J^-)=({Q^\prime}^{1,0})^h \oplus ({Q^\prime}^{0,2})^v\\
\end{align*}

Take a section $f:U\longrightarrow W$ , where $U$ is an open set of $M$ with associated tensors $J^Q\in \Gamma(U,End Q)$ and $F=i(J^Q)\in \Gamma(U,\wedge^2Q)$, the map $i$ is defined as (\ref{i}). Then $f$ is called transversely $J^\pm$- holomorphic if $f:(U,J^\pm)\longrightarrow (W,J^\pm)$ is holomorphic, i.e., $f_*\circ J^Q=J^\pm \circ f_*$.  It can be proved (cf. Proposition 3.2 of \cite{sa}) that a section $f$ is transversely  $J^+$-holomorphic (respectively $J^-$-holomorphic) if and only if $D_1F=0$ (respectively $D_2F=0$).

It is a well-known fact that $J^-$ is never integrable.  Now suppose that the transverse almost complex structure $J^+$ on the normal bundle of the canonical foliation of $W$ is integrable, hence, every point of $W$ belongs to a transversely  $J^+$- holomorphic section $f$ discussed earlier and $D_1F=0$. Let $X,Y\in \Gamma(U,Q^{1,0})$, using Lemma \ref{Lemma 1.2} we have that the bundle $Q^{1,0}$ is preserved by $\nabla_X$, moreover, $[X,Y]\in Q^{1,0}$. Let us denote the transversal Riemannian curvature of $M$ by $R^Q$, then 
\begin{align*}
    R^Q(X,Y)Z=\left([\nabla_X,\nabla_Y]-\nabla_{[X,Y]}\right)Z\in Q^{1,0}.
\end{align*}
where $X,Y,Z\in \Gamma(U,Q^{1,0})$. Hence, $R^Q(X,Y)$ has no component in $Q^{0,2}$, and moreover, $R$ has no components in $Q^{0,2}\otimes Q^{0,2}$, which in the case of $(M,\mathcal{F})$ with foliation of codimension $2q$, for $q\geq 3$ means that $M$ is transversely  conformally flat (see e.g. \cite{N}). Recall that the Hodge star $\star$, is a map from $k$-forms to $(2q-k)$-forms
\begin{equation*}
    \star:\bigwedge\nolimits^k\rightarrow \bigwedge\nolimits^{2q-k}
\end{equation*}
therefore when $q=2$, for $2$-forms we have $\star:\bigwedge^2\rightarrow \bigwedge^2$. In this case, $\star^2=+1$ and it leads to the decomposition of the (transverse) Weyl tensor $\mathcal{W}$ into an anti-self-dual $\mathcal{W}^+$ and a self-dual part $\mathcal{W}^-$.

\begin{definition}
The foliated Riemannian manifold $(M,\mathcal{F})$ is called transversely  self-dual (respectively anti-self-dual) if $\mathcal{W^+}$(respectively $\mathcal{W^-}$) vanishes.
\end{definition}

In the case of non-foliated manifolds the following result was discussed in the paper by M. Atiyah, N.J. Hitchin and I.M. Singer \cite{AHS}, and now using the earlier discussions, we can formulate its counterpart for the case that we have a foliated Riemannian manifold $(M,\mathcal{F})$ with the foliation of codimension $2q$ and its normal twistor bundle, which as before we denote by $W$ and normal twistor bundle has the normal bundle of its foliation denoted by $Q(W)$.

\begin{theorem}\label{th4}
$Q(W,{\mathcal F}^Q,J^+)$ is a tranversely complex foliated manifold if and only if  we have the following
\begin{itemize}
    \item $M$ is transversely  conformally flat for $q\geq3$.
    \item $M$ is transversely  anti-self-dual for $q=2$.
\end{itemize}

\begin{proof}
The proof for the first case (i.e., $q\geq 3$) follows from the earlier discussions, as for when $q=2$, the condition of $R^Q$ having no component in $Q^{0,2}\otimes Q^{0,2}$ means that the self-dual part of the transverse Weyl tensor, which is $\mathcal{W^-}$, vanishes and therefore $M$ is transversely  anti-self-dual. The converse is due to the work done by S. Salamon in \cite{sa}.
\end{proof}
\end{theorem}

\begin{theorem}\label{th5}
Let $(M,g,\mathcal{F},J)$ be transversely  cosymplectic and $(N,\mathcal{F}^\prime)$ a foliated Riemannian manifold with a normal twistor bundle $Z(N,\mathcal{F}_Z^\prime)$ and we have the projection $\pi: Z(N,\mathcal{F}_Z^\prime)\longrightarrow N$. Then if $\psi: M\longrightarrow Z(N,\mathcal{F}^\prime)$ is transversely  $J^-$-holomorphic, $\pi\circ \psi$ is transversely  harmonic.

\begin{proof}
Let us for the sake of brevity denote the normal twistor bundle by $W^\prime$ and the map $\phi:=\pi\circ\psi$. Then the map $\psi$  can be viewed as a section over $M$ of the pull-back bundle $\phi^{-1}W^\prime\subset \phi^{-1}\bigwedge\nolimits^2Q^\prime$, where $Q^\prime=TN/T\mathcal{F}^\prime$ is the normal bundle of the foliation of $N$. Hence, having $\psi$ defines a transverse almost complex structure on the fibers of $\phi^{-1}Q^\prime$, which gives the decomposition 
\begin{equation*}
    (\phi^{-1}Q^\prime)^\mathds{C}={\hat{Q}}^{1,0}\oplus {\hat{Q}}^{0,1}.
\end{equation*}
also the fundamental $2$-vector $F\in\Gamma(M,\phi^{-1}\bigwedge\nolimits^2Q^\prime)$. Moreover, as before we have $\nabla F=D_1F+D_2F$, where
\begin{align}
    D_1F\in (Q^{1,0}\otimes{\hat{Q}}^{2,0}) \oplus (Q^{0,1}\otimes {\hat{Q}}^{0,2}),\\
    D_2F\in  (Q^{1,0}\otimes{\hat{Q}}^{0,2}) \oplus (Q^{0,1}\otimes {\hat{Q}}^{2,0}).
\end{align}
where $Q=TM/T\mathcal{F}$. Let us take $X\in \Gamma Q^{1,0}$, then \begin{align*}
    \psi_*X=(\phi_*X)^h+(\nabla_XF)^h
\end{align*}
which, since $\psi$ is transversely  $J^-$-holomorphic, belongs to $Q^{1,0}(W^\prime,J^-)$, therefore, $\phi_*X\in {\hat{Q}}^{1,0}$ and $\nabla_XF\in{\hat{Q}}^{0,2}$. These results are a consequence of $tr(D_2F^M)=0$ and $D_2F=0$, hence, the rest of the proof follows likewise Theorem $\ref{th 3}$ , having $F$ instead of $F^N$.
\end{proof}
\end{theorem}

Let $(M^{2n+2q},\mathcal{F})$ be a foliated Riemannian manifold with the foliation being transversely  almost Hermitian and having codimension $2q$. For the tangent bundle of $M$ we have the decomposition $TM=L \oplus Q$, where $L=T\mathcal{F}$. Let $J$ be an almost complex structure on $TM$ which preserves the tangent bundle to the leaves, i.e., $J(L)\subset L$. Then  the restriction $J^L$ of $J$ to $L$ defines an almost complex structure in $L.$ Moreover, $J$ induces an almost complex structure $J^Q$ in the normal bundle $Q$. If the almost complex  structure $J$ is compatible with  the bundle-like metric $g$ then $J$  decomposes as 

\begin{align}
    J=J^L\oplus J^Q
\end{align}
as  $J(L)\subset L$ implies that  $J(L^{\perp})\subset L^{\perp}.$ Recall the normal twistor space defined as (\ref{Z}), now using $J^L$ we can define the tangent twistors as

\begin{align}
    Z(\mathcal{F})_x=\{J^L_x:L_x\rightarrow L_x, (J^L)^2=-id_L\}
\end{align}
where $x\in M$ and the bundle of tangent twistors $Z(\mathcal{F})$ as the union of $Z(\mathcal{F})_x$ over all $x\in M$. 

\medskip

On a foliated Riemannian manifold $(M,\mathcal{F})$, there exist three twistor spaces, which we denote by $Z(M)$, $Z(M,\mathcal{F})$ and $Z(\mathcal{F})$ defined respectively using $J$, $J^Q$ and $J^L$. Hence $Z(M)$ is the twistor space of the foliated manifold $(M,\mathcal{F})$, $Z(M,\mathcal{F})$ is the bundle of normal twistors and $Z(\mathcal{F})$ is bundle of tangent twistors, moreover, we have $Z(M)=Z(\mathcal{F})\times_M Z(M,\mathcal{F})$. Similar to the definition of transversely holomorphic map given in Definition \ref{def1} one can define the tangentially holomorphic map.

\begin{definition}
Let $(M,\mathcal{F})$ and $(M^\prime,\mathcal{F}^\prime)$ be two foliated Riemannian manifolds with tangential almost complex structures $J^L$ and $J^{L^{\prime}}$ respectively. The smooth foliated map $\phi:(M,\mathcal{F})\longrightarrow (M^\prime,\mathcal{F}^\prime)$ is tangentially holomorphic if 
\begin{align*}
    \phi^{L}_{*}\circ J^L=J^{L^{\prime}}\circ \phi^{L}_{*}.
\end{align*}
where $\phi^{L}_{*}=\phi_{*}|_L$ , is defined as the restriction of the push-forward map $\phi_{*}$ to the tangents to the foliation  and it is a section of the bundle $L^*\otimes \phi^{-1}L^\prime\longrightarrow M$.
\end{definition}

Let $(M,\mathcal{F})$ be a foliated Riemannian manifold and $Z(M)$ be its twistor space and $\pi:Z(M)\rightarrow (M,\mathcal{F})$ be the projection map. A connection $\omega$ on $SO_g(M)$ induces a splitting of the tangent bundle of the twistor space $TZ(M)=H_M^{(\omega)}\oplus V_M$. At each point $P\in Z(M)$, with $\pi(P)=x\in M$, let $X\in T_PZ(M)$, which splits into horizontal and vertical parts $X=X^h+X^v$. An almost complex structure on $Z(M)$ can be defined as

\begin{align}\label{ACS}
    J_w(X):=(\pi^{-1}_{*|x}\circ P \circ \pi_{*|P})(X^h)+ P \circ X^v
\end{align}

Now suppose we have two Riemannian manifolds $M$ and $M^\prime$ with the same dimensions $2n$ and let $\phi:M\longrightarrow M^\prime$ be a local diffeomorphism, one can lift this map to the corresponding twistor spaces of the manifolds by defining the lifted map as $\Phi(J):=\phi_*\circ J\circ (\phi_*)^{-1}$,as a result we have a diagram

 \begin{displaymath}
 \begin{tikzcd}
   Z(M) \arrow[r, "\Phi "] \arrow[d]
   & Z(M^\prime) \arrow[d]\\
   M \arrow[r, "\phi"]
   & M^\prime
 \end{tikzcd}
 \end{displaymath}

 Recall that there are two almost complex structures on the twistor space by reversing the sign on the vertical part. The map $\Phi$ acts using the map $\phi_*$, hence using (\ref{ACS}), it can be shown that $\Phi$ is holomorphic with respect to the same type of almost complex structure, if and only if the map $\phi$ is holomorphic. Now suppose we have two foliated Riemannian manifolds $(M,\mathcal{F})$ and $(M^\prime,\mathcal{F}^\prime)$ with the same dimensions $2n+2q$ and the transversely Hermitian foliations having codimension $2q$. Recall that on a foliated manifold we have a bundle of normal twistors $Z(M,\mathcal{F})$ and a bundle of tangent twistors $Z(\mathcal{F})$. Using similar argument one can also show that there exists a lifted map for normal twistors (respectively tangent twistors) which is holomorphic if and only if the map between the base manifolds is transversely (respectively tangentially) holomorphic, hence the following diagrams

\begin{figure}
\centering
 \begin{tikzcd}
  Z(M,\mathcal{F}) \arrow[r, "\Phi^Q "] \arrow[d]
   & Z(M^\prime,\mathcal{F}^\prime) \arrow[d]\\
   (M,\mathcal{F}) \arrow[r, "\phi"]
   & (M^\prime,\mathcal{F}^\prime)
 \end{tikzcd}\qquad\begin{tikzcd}
   Z(\mathcal{F}) \arrow[r, "\Phi^L "] \arrow[d]
   & Z(\mathcal{F}^\prime) \arrow[d]\\
   (M,\mathcal{F}) \arrow[r, "\phi"]
   & (M^\prime,\mathcal{F}^\prime)
   \end{tikzcd}{\caption{respectively normal and tangential case. }}
   \end{figure}

 In Theorem \ref{th 3} we proved that having some conditions on the foliated Riemannian manifolds results in any transversely  holomorphic map between them being transversely harmonic, therefore, it is straightforward to imagine a counterpart for this result as following.

\begin{theorem}\label{th6}
Let $(M,\mathcal{F})$ and $(M^\prime,\mathcal{F}^\prime)$ be two transversely almost Hermitian foliated manifolds, which are respectively transversely cosymplectic and transversely $(1,2)$-symplectic. Furthermore, let $Z(M,\mathcal{F)}$ and $Z(M^\prime,\mathcal{F^\prime})$ be the corresponding bundles of normal twistors. Suppose we have the following

\begin{displaymath}
  \begin{tikzcd}
   Z(M,\mathcal{F}) \arrow[r, "\Phi^Q "] \arrow[d]
   & Z(M^\prime,\mathcal{F}^\prime) \arrow[d]\\
   M \arrow[r, "\phi"]
   & M^\prime
 \end{tikzcd}
 \end{displaymath}
 where $\Phi^Q(J^Q):=\phi^T_*\circ J^Q\circ (\phi^T_*)^{-1}$, then if the map $\Phi$ is transversely  holomorphic, the map $\phi$ is transversely  harmonic.

\end{theorem}

It is possible to take the open subsets of leaves of a foliated manifold $(M,\mathcal{F)}$ as the open subsets, hence, a new topology can be introduced and the set $M$ having a compatible differential structure with this new topology is a manifold, which we denote by $M_{\mathcal{F}}$. It can be shown that a smooth map $\phi:(M,\mathcal{F})\longrightarrow (M^\prime,\mathcal{F}^\prime)$ is foliated if and only if it induces a smooth map $\hat{\phi}:M_{\mathcal{F}}\longrightarrow M^\prime_{\mathcal{F}^\prime}$. Let us denote by $\tau(\hat{\phi})$ the tension field of the map $\hat{\phi}$, then $\hat{\phi}$ is called harmonic if $\tau(\hat{\phi})=0$. Moreover, $\phi$ is called leaf-wise harmonic if $\hat{\phi}$ is harmonic.
Now one can also formulate a counterpart for Theorem \ref{th 3} using leaf-wise harmonic and tangentially holomorphic maps. As a final part of this section, we present the following corollary (cf. \cite{kw}, Corollary 2.1), for more details see \cite{kw}.

\begin{corollary}\label{cor3}
Let $\phi:(M,\mathcal{F})\longrightarrow (M^\prime, \mathcal{F}^\prime)$ be a smooth foliated mapping between two regular foliated Riemannian manifolds (i.e. the foliation is regular), with $\mathcal{F}$ being minimal and $\mathcal{F}^\prime$ totally geodesic. Then the map $\phi$ is harmonic if and only if $\phi$ is transversely  harmonic and leaf-wise harmonic.
\end{corollary}

Recall that under some conditions on the base manifolds, the holomorphicity of the maps on the base manifolds of twistor bundles would result in the maps being harmonic. As in Theorem \ref{th6}, the harmonicity of the base map is a result of the holomorphicity of the lifted map on the corresponding twistor bundles which we shall denote by $\Phi^L$ and $\Phi^Q$ for the tangential and normal case respectively. Let us recall some of the known conditions for which a transversely holomorphic map is transversely  harmonic in the following proposition.

\begin{prop}\label{prop2}
Let $(M,\mathcal{F})$ and $(M^\prime,\mathcal{F}^\prime)$ be two almost Hermitian foliated manifolds and $\phi:M\longrightarrow M^\prime$ be a smooth foliated map. Then any transversely holomorphic map $\phi$ is transversely  harmonic if

\begin{enumerate}

    \item $(M,\mathcal{F})$ is transversely cosymplectic and $(M^\prime,\mathcal{F}^\prime)$ is transversely  $(1,2)$-symplectic.
    \item $(M,\mathcal{F})$ and $(M^\prime,\mathcal{F}^\prime)$ are transversely  $(1,2)$-symplectic
    \item $(M,\mathcal{F})$ and $(M^\prime,\mathcal{F}^\prime)$ are Riemannian manifolds with tranversally K\"{a}hler foliations.
\end{enumerate}

\begin{proof}
1. follows from Theorem \ref{th 3}, 2. is an immediate consequence of 1. and for the proof of 3. see \cite{JJ}.
\end{proof}
\end{prop}

 Similarly one can have the counterpart for the tangential case. Corollary \ref{cor3} together with these results can be combined in the following Proposition.
\begin{prop}
Let $\phi:(M,\mathcal{F})\longrightarrow (M^\prime, \mathcal{F}^\prime)$ be a smooth foliated mapping between two regular foliated Riemannian manifolds (i.e. the foliation is regular), with $\mathcal{F}$ being minimal and $\mathcal{F}^\prime$ totally geodesic. Furthermore, let us suppose that we have one of the conditions in the Proposition \ref{prop2} and the tangential counterpart. Then the map $\phi$ is harmonic if and only if the lifted maps $\Phi^L$ and $\Phi^Q$ are respectively tangentially and transversely  holomorphic.
\end{prop}

\section{Orbifolds}\label{S6}

In 1956, I. Satake \cite{sat} introduced a new generalization of the notion of manifolds that he named $V$-manifolds. Currently due to W. Thurston \cite{th}, they  are known as orbifolds and have applications in both mathematics and physics, especially in the string theory. It is a well-known result that the leaf space of a Riemannian foliation with compact leaves is  an orbifold and any orbifold can be realized as the leaf space of a Riemannian foliation, cf. \cite{ha}. We will  reformulate some of the results of the Section \ref{S5} for orbifolds. In this section we follow the notations and borrow some notions from \cite{w2}, which can be consulted for more discussions on the subject.

Let $X$ be a topological space, $\tilde{U}\subset\mathds{R}^n$ be a connected open subset, $\Gamma$ be a finite group of smooth diffeomorphisms of $\tilde{U}$, and $\phi:\tilde{U}\longrightarrow X$ be a map which is $\Gamma$-invariant and induces a homeomorphism of $\tilde{U}/\Gamma$ onto an open subset $U\subset X$. The triple $(\tilde{U},\Gamma,\phi)$ is called an $n$-dimensional orbifold chart on X.

An embedding $\lambda:(\tilde{U},\Gamma,\phi)\longrightarrow (\tilde{V},\Delta,\psi)$ between two orbifold charts is a smooth embedding $\lambda:\tilde{U}\longrightarrow \tilde{V}$ which satisfies $\psi\circ\lambda=\phi$.

Let $\mathcal{A}=\{(\tilde{U}_i,\Gamma_i,\phi_i)\}_{i\in I}$ be a family of such charts, it is called an orbifold atlas on $X$, if it covers $X$ and any two charts are locally compatible in the following sense:
given two charts $\{(\tilde{U}_i,\Gamma_i,\phi_i)\}_{i=1,2}$ and $x\in U_1\cup U_2$, there exists an open neighborhood $U_3\subset U_1\cup U_2$ containing $x$ and a chart $(\tilde{U}_3,\Gamma_3,\phi_3)$, $U_3=\phi_3(\tilde{U}_3)\subset X$ such that it can be embedded into the other two charts. As in the case of manifolds, one can define a maximal atlas.

\begin{definition}
A Hausdorf paracompact topological space $X$ together with a maximal orbifold atlas $\mathcal{A}$ is called a smooth $n$-dimensional orbifold.
\end{definition}

Some of the notions given earlier can be generalized to orbifolds. Let $f:X\longrightarrow Y$ be a map between two Riemannian orbifolds, similar to earlier we have the differential $df:TX\longrightarrow TY$, which is a section of the bundle
\begin{align*}
    TX\otimes f^{-1}TY
\end{align*}
and having the Levi-Civita connections $\nabla^X$ and $\nabla^Y$, one can define a connection $D$ on the bundle $TX\otimes f^{-1}TY$. Therefore we are able to define the tension field $\tau(f)$ for the map between $(X,g_X)$ and $(Y,g_Y)$
\begin{align*}
    \tau(f)=tr(Ddf)
\end{align*}
which is a section of $f^{-1}TY\longrightarrow X$.

\begin{definition}
A complete orbifold mapping $f:X\longrightarrow Y$ between two Riemannian orbifolds $(X,g_X)$ and $(Y,g_Y)$ is harmonic if its tension field vanishes.
\end{definition}

For any orbifold there exists a Riemannian foliation whose leaf space is the orbifold. Suppose $X$ and $Y$ are two orbifolds with $L(X)$ and $L(Y)$ denoting the linear frame bundles respectively. The frame bundles are naturally foliated by the fibers and the foliations are denoted by $\mathcal{F}_X$ and $\mathcal{F}_Y$. Any smooth complete orbifold diffeomorphism $f:X\longrightarrow Y$ can be lifted to a foliated mapping $L(f):L(X)\longrightarrow L(Y)$.

\begin{theorem} \cite{w2}
Let $(X,g)$ and $(Y,h)$ be two Riemannian orbifolds and $f:X\longrightarrow Y$ be a smooth complete mapping. Then $f$ is harmonic if and only if the induced (foliated) mapping $L(f):(L(X),\mathcal{F}_X,g_L)\longrightarrow (L(Y),\mathcal{F}_Y,h_L)$ is transversely  harmonic.
\end{theorem}

The above theorem  shows that as expected some of the properties of the Riemmanian orbifold correspond to the properties of the foliated Riemannian manifold. Now, we will give the counterpart of some the results obtained earlier.

\begin{theorem}
Let $X$ and $Y$ be two almost Hermitian orbifolds, which are respectively cosymplectic and $(1,2)$-symplectic. Then any holomorphic complete orbifold map $\phi:X\longrightarrow Y$ is harmonic.
\end{theorem}

Recall that $(M^{2n},\mathcal{F})$ is a foliated Riemannian manifold with foliation $\mathcal{F}$ and $W:=Z(M,\mathcal{F}_Z)$ is its normal twistor bundle with foliation $\mathcal{F}_Z$. As discussed earlier the leaf spaces $M/\mathcal{F}$ and $W/\mathcal{F}_Z$ correspond to orbifolds, which we denote by $\mathcal{O}_M$ and $\mathcal{O}_W$ respectively. The almost complex structures $J^\pm$, also induce their counterpart on $\mathcal{O}_W$, which by abuse of notation we shall denote them by $J^\pm$. Therefore, we put the following theorem as the orbifold counterpart of Theorem \ref{th4}.

\begin{theorem}
$(\mathcal{O_W},J^+)$ is a complex orbifold if and only if

\begin{itemize}
    \item $\mathcal{O}_M$ is conformally flat for $n\geq 3$.
    \item $\mathcal{O}_M$ is anti-self-dual for $n=2$.
\end{itemize}
\end{theorem}

\begin{prop}
The twistor space of a Riemannian orbifold is a manifold.

\begin{proof}
As we discussed in Section \ref{S3} the twistor space of a manifold $(M^{2n},g)$ is a bundle over $M$ and each point of $Z(M)$ represents an almost complex structure on $M$ which is compatible with the metric and orientation. The dependence of the twistor space on the metric $g$ is in fact conformal, i.e., it depends on the conformal class $[g]$. The counterpart of positive oriented orthonormal frame bundle $SO_g(M)$ can also be constructed on an orbifold, and in fact this bundle is a manifold in this case, therefore, the twistor space of an orbifold which is associated to this bundle is also a manifold. It is worth noting that in a maximal orbifold atlas $\{(\tilde{U}_i,\Gamma_i,\phi_i)\}_{i\in I}$ , if $g_i$ is the Riemannian metric on $\tilde{U}_i$ , for an embedding between two charts we have a map $(\tilde{U}_i,g_i)\rightarrow (\tilde{U}_j,g_j)$, as an isometry which is an even more restrictive condition than being conformal we needed.
\end{proof}
\end{prop}

Taking Theorem \ref{th5} into consideration and denoting the orbifolds obtained from $(N,\mathcal{F}^\prime)$ and $Z(N,\mathcal{F}_Z^\prime)$ by $\mathcal{O}_N$ and $\mathcal{O}_{W^\prime}$, respectively, we get the following theorem.

\begin{theorem}
Let $\mathcal{O}_M$ be a cosymplectic orbifold, having the orbifolds $\mathcal{O}_N$ and $\mathcal{O}_{W^\prime}$ and the projection $\pi: \mathcal{O}_{W^\prime}\longrightarrow \mathcal{O}_N$, if $\psi: \mathcal{O}_M\longrightarrow \mathcal{O}_{W^\prime}$, is $J^-$-holomorphic, $\pi\circ \psi$ is harmonic.
\end{theorem}

\section*{Acknowledgments}
The first  author acknowledges support of the Polish Ministry of Science and
Higher Education.

\vskip0.5cm

\noindent Rouzbeh Mohseni

\noindent  Jagiellonian University in Krakow, Institute of Mathematics, ul. St. Lojasiewicza 4, 30-348 Krakow, Poland.

\noindent {\it E-mail address} : \textbf{rouzbeh.mohseni@doctoral.uj.edu.pl}

\vskip0.3cm

\noindent Robert A. Wolak

\noindent Jagiellonian University in Krakow, Institute of Mathematics, ul. St. Lojasiewicza 4, 30-348 Krakow, Poland.

\noindent{\it E-mail address} : \textbf{robert.wolak@uj.edu.pl}

\end{document}